\documentclass[12pt, cite, epsf]{article}
\usepackage{amsfonts}
\usepackage{mathrsfs}

\usepackage{amstext, graphicx, amsmath, latexsym, amssymb, amscd, amsfonts}

\newtheorem{theorem}{Theorem}[section]
\newtheorem{corollary}{Corollary}[section]
\newtheorem{lemma}{Lemma}[section]
\newtheorem{proposition}{Proposition}[section]
\newtheorem{remark}{Remark}[section]
\newtheorem{example}{Example}[section]

\newcommand{\R}{{\mathbb{R}}}

\DeclareMathOperator*\argmax{arg\,max}
\DeclareMathOperator*\co{conv}
\DeclareMathOperator*\cl{cl}
\DeclareMathOperator*\pos{pos}
\DeclareMathOperator*\inte{int}
\DeclareMathOperator*\ri{ri}

\DeclareMathOperator*\es{end}
\DeclareMathOperator*\ebm{ebm}
\DeclareMathOperator*\bd{bdry}

\begin{document}

\begin{center}

%
%\begin{center}
%\today
%\end{center}

\begin{center}
{\Large\Large\sc {\bf On Error Bound Moduli for Locally Lipschitz and Regular Functions}}
\end{center}

{\sc M.\ H.\ Li}
\\ {\small School of Mathematics and Finance, Chongqing University of Arts and Sciences, Yongchuan, Chongqing, 402160, China}\\
Email: minghuali20021848@163.com\\[0.5cm]

{\sc K.\ W.\ Meng}
\\ {\small School of Economics and Management, Southwest Jiaotong University, Chengdu 610031, China}\\
Email: mkwfly@126.com\\[0.5cm]

{\sc X. \ Q.\ Yang}\\ {\small Department of Applied Mathematics, The Hong Kong Polytechnic University,  Hong Kong}\\
Email: mayangxq@polyu.edu.hk
\end{center}

 \noindent{\bf Abstract}:  In this paper we study local error bound moduli for a locally Lipschitz and regular function  via its outer limiting subdifferential set.  We show that  the distance of 0 from the outer limiting subdifferential of the support function of the subdifferential set, which is essentially the distance of 0 from
 the end set of the subdifferential set, is an upper estimate of the local error bound modulus. This upper estimate becomes tight for a convex function under some regularity conditions.  We show that the distance of 0 from the outer limiting subdifferential set of a  lower $\mathcal{C}^1$ function   is equal to the local error bound modulus.


\section{Introduction}

Error bounds play a key role in variational analysis. They are of great importance for subdifferential calculus, stability and sensitivity analysis, exact penalty functions, optimality conditions, and convergence of numerical methods, see the excellent survey papers \cite{ay,lew,pang} for more details. It should be noticed that the notion of error bounds is closely related to some other important concepts: weak sharp minima, calmness and metric subregularity, see \cite{ay1,bur0,do,io,kla,ng,roc}.

In this paper, we study  local error bound moduli  in finite dimensional spaces.
 We say that  a function $\phi:\R^n\to \overline{\R}:=\R\cup\{+ \infty\}$ has a local error bound at $\bar{x}\in [\phi\leq 0]$ if there exist some $\tau>0$ and some neighborhood $U$ of $\bar{x}$ such that
\begin{equation}\label{eq1}
\tau d(x,[\phi\leq 0])\leq \phi(x)_+\;\forall x\in U,
\end{equation}
where $[\phi\leq 0]:=\{x\in \R^n|\phi(x)\leq0\}$ and $t_+:=\max\{t,0\}$ for all $t\in \R$. The supremum of all possible constants $\tau$ in (\ref{eq1}) (for some associated $U$) is called the local error bound modulus of $\phi$ at $\bar{x}$, denoted by $\ebm(\phi, \bar{x})$. We define $\ebm(\phi, \bar{x})$  as $0$ if $\phi$ does not have a local error bound at $\bar{x}$. Clearly, the local error bound modulus of $\phi$ at $\bar{x}$ can be alternatively defined as
\[
 \ebm(\phi,\bar{x})=\liminf_{x\rightarrow \bar{x}, \phi(x)>0}\frac{\phi(x)}{d(x, [\phi\leq 0])}.
\]
We know from the literature \cite{knt10,fhko10,io15} that the distance of 0 from the outer limiting subdifferential of a lower semicontinuous (lsc) function $\phi$ at $\bar{x}$, is a lower estimate of $\ebm(\phi, \bar{x})$, which becomes tight when $\phi$ is convex.

In this paper we consider the local error bound modulus of a locally Lipschitz and regular function $\phi$ and
 establish that the distance of 0 from the outer limiting subdifferential of the support function of the subdifferential
 $\partial \phi(\bar{x})$ at 0 is an upper estimate of $\ebm(\phi, \bar{x})$. We also investigate the geometric structure of this outer limiting subdifferential and show that it is equal to the closure of the end set of the subdifferential $\partial \phi(\bar{x})$, while the closure is surplus when the subdifferential set is a polyhedron. Thus the upper estimate is essentially the distance of 0 from the end set of $\partial \phi(\bar{x})$.
We prove that, for convex function $\phi$, under ACQ and ETA (see Remark \ref{rem-acq+eta}), the upper estimate is tight. To the best of our knowledge, the first result of the kind is that \cite{hu2} proved that for sublinear function $\phi$, the $\ebm(\phi, \bar{x})$ is equal to the distance of 0 from the end set of $\partial \phi(\bar{x})$.

For  lower $\mathcal{C}^1$ function $\phi$, we show that the distance of 0 from the outer limiting subdifferential of $\phi$ at $\bar{x}$ is equal to $\ebm(\phi, \bar{x})$. This generalizes the corresponding results in  \cite{knt10,fhko10,io15} for convex function $\phi$.

Throughout the paper we use the standard notations of variational analysis; see the seminal book \cite{roc} by Rockafellar and Wets. Let $A\subset \R^n$. We denote the closure, the boundary, the convex hull and the positive hull of $A$ respectively by $\cl A$, $\bd A$, $\co A$ and $\textup{pos}A:=\{0\}\cup \{\lambda x| x\in A$ and $\lambda>0\}$.

The Euclidean norm of a vector $x$ is denoted by $||x||$, and the inner product of vectors $x$ and $y$ is denoted by $\langle x, y\rangle$. Let $B(x,\varepsilon)$ denote a neighborhood of $x\in \R^n$ with the radius $\varepsilon >0$.
We say that $A$ is locally closed at a point $x\in A$ if $A\cap U$ is closed for some closed neighborhood $U$ of $x$. The polar cone of $A$ is defined by
\[
A^*:=\{v\in \R^n| \langle v, x\rangle\leq 0\;\forall x\in A\}.
\]
The support function $\sigma_A: \R^n\to \bar{\R}$ of $A$ is defined by
\[
\sigma_A(w):=\sup_{x\in A}\langle x, w\rangle.
\]
For a closed and convex set $A$ with $0\in A$, the gauge of $A$ is the function $\gamma_A:\R^n\to \overline{\R}$ defined by
\[
\gamma_A(x):=\inf\{\lambda\geq 0| x\in \lambda A\}.
\]

The distance of $x$ from $A$ is defined by
\[
d(x,A):=\inf_{y\in A}||y-x||.
\]
For $A=\emptyset,$ we define $d(x,A)=+\infty$. The projection mapping $P_A$ is defined by
\[
P_A(x):=\{y\in A|||y-x||=d(x,A)\}.
\]
Let $x\in A$. We use $T_A(x)$ to denote the tangent cone to $A$ at $x$, i.e. $w\in T_A(x)$ if there exist sequences $t_k\downarrow 0$ and $\{w_k\}\subset \R^n$ with $w_k\rightarrow w$ and $x+ t_k w_k\in A \;\forall k$. We denote by $N_A^P(x)$ the proximal normal cone to $A$ at $x$, i.e., $v\in N_A^P(x)$ if there exists some $t>0$ such that $x\in P_A(x+tv)$. The regular normal cone $\hat{N}_A(x)$ to $A$ at $x$ is the polar cone of $T_A(x)$. A vector $v\in \R^n$ belongs to the normal cone $N_A(x)$ to $A$ at $x$, if there exist sequences $x_k\to x$ and $v_k\to v$ with $x_k\in A$ and $v_k\in \hat{N}_A(x_k)$ for all $k$. The set $A$ is said to be regular at $x$ in the sense of Clarke if it is locally closed at $x$ and $\hat{N}_A(x)=N_A(x)$.

Let $g:\R^n\to \bar{\R}$ be an extended real-valued function and $x$ a point with $g(x)$ finite. We denote by $\ker g:=\{x\in \R^n| g(x)=0\}$ the kernel of $g$. The epigraph of $g$ is the set
\[
\textup{epi}g:=\{(x,\alpha)\in \R^n\times \R|g(x)\leq \alpha\}.
\]
It is well known that $g$ is lower semicontinuous (lsc) if and only if $\textup{epi} g$ is closed.
The vector $v\in \R^n$ is a regular subgradient of $g$ at $x$, written $v\in \hat{\partial} g(x)$, if
\[
g(y)\geq g(x)+\langle v, y-x\rangle+o(||y-x||).
\]
The vector $v\in \R^n$ is a (general) subgradient of $g$ at $x$, written $v\in \partial g(x)$, if there exist sequences $x_k\to x$ and $v_k\to v$ with $g(x_k)\to g(x)$ and $v_k\in \hat{\partial} g(x_k)$. The outer limiting subdifferential of $g$ at $\bar{x}$ is defined in  \cite{ knt10, fhko10, io15}  by
\[
\partial^>g(\bar{x})=\{\lim_{k\rightarrow +\infty}v_k\mid \exists x_k\rightarrow_g \bar{x},\,g(x_k)>g(\bar{x}),\,v_k\in\partial g(x_k)\}.
\]

The subderivative function $dg(x):\R^n\to \bar{\R}$ is defined by
\[
dg(x)(w):=\liminf_{t\downarrow 0, w'\to w}\frac{g(x+tw')-g(x)}{t}\;\forall w\in \R^n.
\]
Note that the subderivative $dg(x)$ is a lsc and positively homogeneous function and that the regular subdifferential set can be derived from the
subderivative as follows:
\[
\hat{\partial} g(x)=\{v\in \R^n| \langle v, w\rangle\leq dg(x)(w)\;\forall w\in \R^n\}.
\]
The function $g$ is said to be (subdifferentially) regular at $x\in \R^n$ if $\textup{epi}g$ is regular in the sense of Clarke at $(x,g(x))$ as a subset of $\R^n\times \R$.

For a sequence $\{A_k\}$ of subsets of $\R^n$, the outer limit $\limsup_{k\to \infty}A_k$ is the set consisting of all possible cluster points of sequences $x_k$
with $x_k\in A_k$ for all $k$, whereas the inner limit $\liminf_{k\to \infty}A_k$ is the set consisting of all possible limit points of such sequences. $\{A_k\}$ is said to converge to $A\subset \R^n$ in the sense of Painlev\'{e}-Kuratowski, written $A_k\to A$, if
\[
\limsup_{k\to \infty}A_k=\liminf_{k\to \infty}A_k=A.
\]

For a set-valued mapping $S:\R^n\rightrightarrows \R^m$ and a point $\bar{x}\in \R^n$, the outer limit of $S$ at $\bar{x}$ is defined by
\[
\limsup_{x\to \bar{x}}S(x):=\{u\in \R^m| \exists x_k\to \bar{x}, \exists u_k\to u\;\textup{with}\;u_k\in S(x_k)\}.
\]
$S$ is outer semicontinuous (osc, for short) at $\bar{x}$ if and only if
\[
\limsup_{x\to \bar{x}}S(x)\subset S(\bar{x}).
\]

A face of a convex set $A$ is a convex subset $A'$ of $A$ such that every closed line segment in $A$ with a relative interior point in $A'$ has both endpoints in $A'$.  An exposed face of $A$ is the intersection of $A$ and a non-trivial supporting hyperplane to $A$. See \cite{rock70}. For a nonempty and convex set $A\subset \R^n$, the end set of $A$ is defined in \cite{hu, hu1} by
\[
\es(A):=\{x\in \textup{cl}A| tx\not\in \textup{cl}A\;\forall t>1\}.
\]

\section{Support Function and Its Outer Limiting Subdifferential}

In this section, we  study the outer limiting subdifferential of the support function of a general compact and convex set and show that it is the
 closure of its end set. In the next section, we shall apply these results  to study the error bound modulus for a locally Lipschitz and regular function
      as its  subdifferential set is compact and convex.
It is worth noting that the results presented in this section
 have their own interest in the field of convex analysis and optimization.

Let $C$ be a compact and convex subset of $\R^n$.  To begin with, we note from \cite[Theorem 8.24, Proposition 8.29, Corollary 8.25]{roc} that
\begin{equation}\label{ccc}
C=\{v\in \R^n\mid \langle v, w\rangle\leq \sigma_C(w)\;\forall w\}
\end{equation}
and
 \begin{equation}\label{supp_part}
 \partial \sigma_C(w)=\argmax_{v\in C}\langle v, w\rangle=C\cap\{v\in \R^n\mid  \langle v, w\rangle=\sigma_C(w)\}.
 \end{equation}
Since $C$ is compact and convex, we have $ \sigma_C(w)<+\infty$ and $ \partial \sigma_C(w)\not=\emptyset$ for all $w\in\R^n$.

We first consider the case that $0 \in C$ and then the general case that $C$ may not contain $0$. Some basic properties of the end set of $C$ are listed in the following lemma.

\begin{lemma}
If $C\subset \R^n$ is compact and convex with $0\in C$, then the following properties hold:
 \begin{description}
   \item[(i)]  $\es(C)\cap \ri C=\emptyset$;
   \item[(ii)] $C=\bigcup_{v\in \es(C)}[0, v]$;
   \item[(iii)] For a subset $ E\subset C$,   $C=\bigcup_{v\in E}[0, v]$  if and only if $\es(C)\subset E$;

   \item[(iv)] $F$ is a nonempty exposed face of $C$ if and only if $F=\partial \sigma_C(w)$ for some $w\not=0$.
 \end{description}
\end{lemma}
\noindent{\bf Proof.} (i) Suppose by contradiction that $v\in \es(C)\cap \ri C$. By the relative interior criterion \cite[Exercise 2.41]{roc}, there must exist some $v'\in C$ such that $v\in \ri[0, v']$, which contradicts to the fact that $v\in \es(C)$.

(ii) This  equality holds because $C$ is convex and compact with $0\in C$.

(iii) The `if' part is trivial due to (ii). As for the `only if' part, we only need to show $\es (C)\subset E$. Let $v\in \es(C)$. Since $C$ is compact, it is clear that $v\in C$. By $C=\bigcup_{v\in E}[0, v]$, there exists some $v'\in E$ such that $v\in [0, v']$. By the definition of the end set, we have $v=v'$. This entails that $\es(C)\subset E$.

(iv) Clearly, any  $\partial \sigma_C(w)$ with $w\not=0$ is an exposed face of $C$. Conversely, if $F\not=\emptyset$ is  exposed in $C$,  then by definition  there exist some  $w\not=0$ and   $\alpha\in \R$ such that
$F=C\cap\{v\in \R^n \mid \langle w, v\rangle=\alpha\}$ and $C\subset \{v\in \R^n\mid \langle w, v\rangle\leq \alpha\}$. The latter inclusion holds if and only if $\sigma_C(w)\leq \alpha$. In view of (\ref{ccc}) and the fact that $F\subset C$, we have $\alpha=\langle w, v\rangle\leq \sigma_C(w)$ for each $v\in F$. This entails that $\alpha=\sigma_C(w)$. In view of (\ref{supp_part}), we have $F=\partial \sigma_C(w)$. The proof is completed. \hfill $\Box$

Throughout this section, we use the following notation:
 \[
 S:=\displaystyle\bigcup_{\sigma_C(w)>0}\partial \sigma_C(w).
 \]
According to Lemma 2.1 (iv),  $S$ is the union of all the exposed faces $\partial \sigma_C(w)$ of $C$ with $\sigma_C(w)>0$.

In next lemma, we prove that $S$ lies in the $\es(C)$, both sharing with the same closure that happens to be the outer limiting subdifferential $\partial^>\sigma_C(0)$.

\begin{lemma}\label{theo-end-set-genel}
If $C\subset \R^n$ is compact and convex with $0\in C$, then
\[
S\subset \es(C)=\gamma^{-1}_C(1)\subset \partial^>\sigma_C(0)=\cl S,
 \]
 entailing that
 \[
 \cl(\es(C))=\cl(\gamma_C^{-1}(1))=\partial^>\sigma_C(0).
 \]
\end{lemma}
\noindent{\bf Proof.}  First, we show   $\es(C)=\gamma^{-1}_C(1)$ and $\cl S=\partial^>\sigma_C(0)$. The first equality follows from the definitions of the end set and the gauge function,  while the second equality  follows readily from  the positive homogeneity of $\sigma_C$ and the definition of outer limiting subdifferential.

  Next, we show $S\subset \es (C)$.   Let $v\in S$, i.e., $v\in  \partial \sigma_C(w)$ for some $w\in \R^n$ with $\sigma_C(w)>0$. In view of (\ref{supp_part}), we have $v\in C$ and $\langle v, w\rangle=\sigma_C(w)$, implying that $\langle tv, w\rangle>\sigma_C(w)$ for all $t>1$. By (\ref{ccc}),  we have  $tv\not\in C$ for all $t>1$. That is,  $v\in \es(C)$.

 Finally we show $\es(C)\subset \cl S$. To begin with, we show that $\es(C)\subset S\cup S^0$, where
   \begin{equation}\label{s0}
   S^0=\bigcup_{\sigma_C(w)=0,\;\partial \sigma_C(w)\not=C} \partial \sigma_C(w).
   \end{equation}
            Let $v\in \es(C)$.  By Lemma 2.1 (i), $v\not\in \ri C$. It then follows from \cite[Theorem 11.6]{rock70} that there exists a non-trivial supporting hyperplane $H$ to $C$ containing $v$. That is, we can find an exposed face
             $F:=C\cap H$ of $C$ such that $v\in F$ and $F\not=C$. By Lemma 2.1 (iv), we can find some $w\not=0$ such that $F=\partial \sigma_C(w)$. This entails that $v\in S\cup S^0$. Therefore, we have  $\es(C)\subset S\cup S^0$ as expected. By Lemma 2.1 (iii), we have  $C=\bigcup_{v\in S\cup S^0}[0, v]$. Observing  that
            $0\in \partial \sigma_C(w)$ for all $w\in \R^n$ with $\sigma_C(w)=0$, we have $\bigcup_{v\in S^0}[0, v]=S^0$. This entails that $C=A\cup S^0$, where  $A:=\bigcup_{v\in S}[0, v]$. Since each $\partial \sigma_C(w)$ with $\sigma_C(w)=0$ and $\partial \sigma_C(w)\not=C$
            is a non-trivial exposed face of $C$, we confirm that $\ri C\cap S^0=\emptyset$ (implying that $\ri C\subset A$ and hence $C\subset \cl A$).  Clearly, we have $A\subset C$ and hence $\cl A\subset C$. That is, we  have $C=\cl A$. On the basis of the fact that $S\subset C$ is bounded,  it's  easy to verify that $\cl A=\bigcup_{v\in \cl S}[0, v]$. Thus, we have $C=\bigcup_{v\in \cl S}[0, v]$. By Lemma 2.1 (iii) again, we   have  $\es(C)\subset \cl S$.

           To sum up, we have shown $S\subset \es(C)=\gamma^{-1}_C(1)\subset \partial^>\sigma_C(0)=\cl S$, which clearly implies that $\cl(\es(C))=\cl(\gamma_C^{-1}(1))=\partial^>\sigma_C(0)$.
           The proof is completed. \hfill $\Box$

\begin{remark}
The closure operation in the equality $\partial^>\sigma_C(0)=\cl (\es(C))$ cannot in general be dropped, because $\partial^>\sigma_C(0)$ is always closed but $\es(C)$ may not be closed, taking for example the simple set $C=\{x\in \R^2\mid 0\leq x_1\leq 1, x_1^2\leq x_2\leq x_1\}$.
\end{remark}

Under some further conditions on the faces of $C$, we show that  $S$ coincides with $\es(C)$.

\begin{lemma}\label{theo-end-set-S}
Assume that $C\subset \R^n$ is compact and convex with $0\in C$.  If, for any $w\in \R^n$ with $\sigma_C(w)=0$ and $\partial \sigma_C(w)\not=C$, all the faces of $\partial \sigma_C(w)$  containing no 0 are exposed in $C$,  then
\[
\es(C)=S.
  \]
  In particular, if $C$ is a polyhedral set, then
  \[
  S=\es(C)=\gamma^{-1}_C(1)=\partial^>\sigma_C(0)=\cl S,\]
   implying that the sets $S$, $\es(C)$ and $\gamma^{-1}_C(1)$ are all closed.
   \end{lemma}
\noindent{\bf Proof.}     We first  show the equality   $\es(C)=S$ under the assumed conditions on the faces of $C$.  Let $v\in \es(C)\cap S^0$,  where $S^0$ is given by (\ref{s0}).   By the definition of $S^0$,  there exists some $w\in \R^n$ with $\sigma_C(w)=0$ and $\partial \sigma_C(w)\not=C$  such that $v\in \partial \sigma_C(w)$.    By \cite[Theorem 18.2]{rock70}, there exists a unique face $F$ of $C$ such that $v\in \ri F$.
      By \cite[Theorem 18.1]{rock70}, we have $F\subset \partial \sigma_C(w)\subset C$.  Clearly, $F$ is also a face of $\partial \sigma_C(w)$.
      We claim that $0\not\in F$,  for otherwise there must exist some $v'\in F$ such that $v\in \ri[0, v']$ (so that $t_0v\in  F\subset C$ for some $t_0>1$),
       contradicting to the assumption that $v\in \es(C)$.
That is,  $F$ is a face of $\partial \sigma_C(w)$ containing no 0, which  is assumed to be exposed in $C$. It then follows from Lemma 2.1 (iv) that $F=\partial \sigma_C(w')$ for  some $w'\not=0$.  As $0\not\in F$, we have $\sigma_C(w')>0$, implying  that $F\subset S$.  Then, we have $v\in S$. This entails that $\es(C)\cap S^0\subset S$. As we have shown in the proof of Lemma 2.2 that $\es(C)\subset S\cup S^0$ and $S\subset \es(C)$, we get the equality  $\es(C)=S$.

To complete the proof, it suffices to note that any polyhedral set has only finitely many faces and all non-trivial faces are exposed ones.
 The proof is completed. \hfill $\Box$

 \begin{remark}
 Without the conditions imposed on faces of $C$ in Lemma \ref{theo-end-set-S}, the union set $S$ may not  be closed as can be seen from Example \ref{exam-S-not-closed} below, demonstrating that the closure operation in the equality $\partial^>\sigma_C(0)=\cl S$  cannot in general be dropped, and that the equality $\es(C)=S$ does not hold in general.
\end{remark}

\begin{example}\label{exam-S-not-closed}  Let $C=\{x\in \R^2\mid 0\leq x_1\leq 2, 0\leq x_2\leq 1+\sqrt{1-(x_1-1)^2}\}$. Clearly, $C$ is a compact and convex set with $0\in C$. By some direct calculation, we have
\[
\es(C)=\{x\in \R^2\mid 0\leq x_1\leq 2, x_2=1+\sqrt{1-(x_1-1)^2}\}\cup \{x\in \R^2\mid x_1=2, 0\leq x_2<1\},
 \]
 while
 \[
 S=\{x\in \R^2\mid 0<x_1\leq 2, x_2=1+\sqrt{1-(x_1-1)^2}\}\cup \{x\in \R^2\mid x_1=2, 0\leq x_2<1\}.
 \]
Let $v=(0, 1)^T$ and   $w=(-1, 0)^T$.  It is easy to verify that $v\in \es(C)\backslash S$, and that the singleton set  $\{v\}$ is a face of $\partial \sigma_C(w)=\{v\in \R^2\mid v_1=0, 0\leq v_2\leq 1\}$ with $\sigma_C(w)=0$, but it is not exposed in $C$.
\end{example}

In the following lemma, we present some equivalent conditions for the equality $\es(C)=\partial^>\sigma_C(0)$ to hold.

\begin{lemma}\label{theo-end-set-closed}
 If $C\subset \R^n$ is compact and convex with $0\in C$, then the following properties are equivalent:
     \begin{description}
       \item[(i)] $\es(C)=\partial^>\sigma_C(0)$;
       \item[(ii)] $\es(C)$ is closed;
       \item[(iii)] $\gamma_C$ is continuous at every $x\in \pos C$ relative to $\pos C$;
       \item[(iv)] $C$ is a radiative subset of $\pos C$ in the sense of \cite[Definition 4.1]{r00}.
     \end{description}
     \end{lemma}
\noindent{\bf Proof.}  The equivalence of (i) and (ii) follows directly from Lemma 2.2, while the equivalence of (iii) and (iv) can be found in \cite[Proposition 4.2]{r00}. It remains to show the equivalence of (ii) and (iii).

[(ii)$\Longrightarrow$(iii)]: Since $0\not\in \es(C)$ and $\es(C)$ is closed, we have $d(0, \es(C))>0$. From \cite[Theorem 4.1]{meng1}, it then follows  that $\pos C$ is closed,   $\gamma_C$ is continuous at 0 relative to $\pos C$, and there is no convergent sequence $\{x_k\}\subset \pos C$ such that $\gamma_C(x_k)\rightarrow +\infty$. Let $x\in \pos C$ with $x\not=0$ and let $x_k\rightarrow x$ with $x_k\in \pos C$ for all $k$. Without loss of generality, we may assume that $x_k\not=0$ for all $k$ and that $\gamma_C(x_k)\rightarrow \beta$ (Note that the sequence $\gamma_C(x_k)$ is bounded). As $\gamma_C$ is lower semi-continuous, we have $\beta\geq \gamma_C(x)>0$.    Since  $\gamma_C(x_k/\gamma_C(x_k))=1$,  we have $x_k/\gamma_C(x_k)\in \es(C)$.  Since $\es(C)$ is assumed to be closed, we have $x_k/\gamma_C(x_k)\rightarrow x/\beta\in \es(C)$. Thus, we have $\gamma_C(x/\beta)=1$ or $\beta=\gamma_C(x)$. This entails that $\gamma_C$ is continuous at $x$ relative to $\pos C$. Therefore, we have (ii)$\Longrightarrow$(iii).

[(iii)$\Longrightarrow$(ii)]: From \cite[Theorem 4.1]{meng1}, it follows that $\pos C$ is closed and $d(0,\es(C))>0$. Let $v_k\rightarrow v$ with $v_k\in \es(C)$ (that is, $\gamma_C(v_k)=1$) for all $k$. It's easy to verify that $v_k\in \pos C$ for all $k$ and  $v\in \pos C$. Moreover, we have $v\not=0$, for otherwise we  have $d(0,\es(C))=0$, a contradiction. By (iii), we have $\gamma_C(v_k)\rightarrow \gamma_C(v)$, implying that $\gamma_C(v)=1$ or equivalently $v\in \es(C)$. This entails the closedness of $\es(C)$. The proof is completed. \hfill $\Box$

Now we present the results similar to the ones in Lemmas 2.2-2.4, for the case when $C$ may not contain $0\in \R^n$.

\begin{theorem}\label{theo-end-set-generel}
Let $C\subset \R^n$ be a compact and convex set not necessarily containing 0, and let $C':=\co(C\cup\{0\})$. Then the following properties hold:
\begin{description}
  \item[(a)] $S\subset \es(C)=\gamma^{-1}_{C'}(1)\subset \partial^>\sigma_C(0)=\cl S$, entailing that  $\cl(\es(C))=\cl(\gamma_{C'}^{-1}(1))=\partial^>\sigma_C(0)$.
      \item[(b)] If, for any $w\in \R^n$ with $\sigma_C(w)=0$ and $\partial \sigma_{C'}(w)\not=C'$, all the faces of $\partial \sigma_{C'}(w)$ containing no 0 are exposed in $C'$,  then $\es(C)=S$.  In particular, if $C$ is a polyhedral set, then $S=\es(C)=\gamma^{-1}_{C'}(1)=\partial^>\sigma_C(0)=\cl S$, implying that the sets $S$, $\es(C)$ and $\gamma^{-1}_{C'}(1)$ are all closed.
    \item[(c)] The following properties are equivalent:
     \begin{description}
       \item[(c1)] $\es(C)=\partial^>\sigma_C(0)$;
       \item[(c2)] $\es(C)$ is closed;
      \item[(c3)] $\gamma_{C'}$ is continuous at every $x\in \pos C'$ relative to $\pos C'$;
       \item[(c4)] $C'$ is a radiative subset of $\pos C'$ in the sense of \cite[Definition 4.1]{r00}.
     \end{description}
   \end{description}
\end{theorem}
\noindent{\bf Proof.} Clearly, $C'=\cup_{0\leq \lambda\leq 1}\lambda C$ is a  compact and convex set with $0\in C'$, and $C'$ is  polyhedral  if   $C$ is   polyhedral. Moreover, it is easy to verify that $\sigma_{C'}(w)=\max\{\sigma_{C}(w), 0\}$ for all $w\in \R^n$,  and that
 \[
 \partial \sigma_{C'}(w)=\left\{
 \begin{array}{ll}
 \partial \sigma_{C}(w)&\mbox{if}\;\sigma_C(w)>0,\\
  \cup_{0\leq \lambda\leq 1}\lambda\partial \sigma_{C}(w)&\mbox{if}\;\sigma_C(w)=0,\\
  \{0\}&\mbox{if}\;\sigma_C(w)<0.\\
 \end{array}
 \right.
 \]
This entails that  $S=\cup_{\sigma_C(w)>0}\partial \sigma_C(w)=\cup_{\sigma_{C'}(w)>0}\partial \sigma_{C'}(w)$.
By definition, we have $\es(C')=\es(C)$.  All results then  follow readily from Lemmas 2.2-2.4.  \hfill $\Box$

By applying Theorem  \ref{theo-end-set-generel}, we can give some formulas for calculating $\partial^>\sigma_C(0)$ when $C$ is the convex hull of a compact subset of $\R^n$.
\begin{corollary}\label{conv-hull-compa}
Let $A$ be a nonempty compact subset of $\R^n$ such that $C=\co A$. In terms of $\mathcal{A}=\{\argmax_{a\in A}\langle a, w\rangle\mid \max_{a\in A}\langle a, w\rangle>0\}$, we have
\begin{equation}\label{tenggeer}
\displaystyle\bigcup_{A'\in \mathcal{A}}\co A'\subset \es(C)=\gamma_C^{-1}(1)\subset \cl\left(\displaystyle\bigcup_{A'\in \mathcal{A}}\co A'\right)=\partial^>\sigma_C(0).
\end{equation}
If  $A$ is a finite set, all the inclusions in (\ref{tenggeer}) become equalities.
\end{corollary}
\noindent{\bf Proof.} It suffices to show that  $A'\in \mathcal{A}$ if and only if there is some $w\in \R^n$ such that $\sigma_C(w)>0$ and $\co A'=\partial \sigma_C(w)$, and then apply Theorem \ref{theo-end-set-generel} in a straightforward way.  \hfill $\Box$

\begin{remark}\label{rem-canovas}
It is easy to verify that $\mathcal{A}$ can be rewritten as
\[
\{A'\subset A \mid \exists w\in \R^n\,\mbox{such that}\, \langle a, w\rangle=1\,\forall a\in A',\; \langle a, w \rangle<1\,\forall a\in A\backslash A'\},
\]
which is in the spirit of the index collection defined   in  C\'{a}novas et al. \cite{ca2} for the case that $A$ is a finite set. On the other hand, when $A$ is a finite set, the equalities in (\ref{tenggeer}) provide a complete characterization of the set $\es(C)$. From which, it is easy to see that
\begin{equation}\label{zhengnh} d(0, \es(C)) > 0.\end{equation}
It is worth noting that (\ref{zhengnh}) has been proved in \cite{zheng,hu}.
\end{remark}

\section{Main Results}

Throughout this section, for a given lsc function  $\phi:\R^n\rightarrow \overline{\R}$, which  is  regular and locally Lipschitz continuous at $\bar{x}$, a point on the boundary of the level set $[\phi\leq 0]$, we shall conduct some variational analysis on $\ebm(\phi, \bar{x})$, the error bound modulus  of $\phi$ at $\bar{x}$. We first show that the distance of 0 from $\partial^>\phi(\bar{x})$,  the outer limiting subdifferential of $\phi$ at $\bar{x}$,  is a lower estimate of $\ebm(\phi, \bar{x})$, while the distance of 0 from $\partial^>\sigma_{\partial\phi(\bar{x})}(0)$,  the outer limiting subdifferential of  $\sigma_{\partial \phi(\bar{x})}$ (the support function of $\partial \phi(\bar{x})$)  at 0,  is an upper estimate of $\ebm(\phi, \bar{x})$. We then show that the lower estimate is tight for a  lower $\mathcal{C}^1$ function  and the upper estimate is tight for a convex function under some regularity conditions.

To begin with, we recall that the inequality
\begin{equation}\label{lower-bound-gen}
\ebm(f, x)\geq d\left(0, \partial^>f(x)\right)
\end{equation}
holds for a lsc function $f$ on $\R^n$ and a point $x$ with $f(x)$ finite, and   the equality
\begin{equation}\label{sharp-bound-convex}
\ebm(f, x)=d\left(0, \partial^>f(x)\right)
\end{equation}
holds if, in addition, $f$ is convex. See \cite{knt10, fhko10,io15}.

\begin{theorem}\label{theo-upper-bound}
Consider a function $\phi:\R^n\rightarrow \overline{\R}$ and a point $\bar{x}$ on the boundary of the level set $[\phi\leq 0]$. If $\phi$ is  regular and locally Lipschitz continuous at  $\bar{x}$, then
\begin{equation}\label{lower-upper-esti}
d\left(0, \partial^>\phi(\bar{x})\right)\leq \ebm(\phi, \bar{x}) \leq d(0,  \partial^>\sigma_{\partial \phi(\bar{x})}(0)).
\end{equation}
\end{theorem}
\noindent{\bf Proof.} In view of (\ref{lower-bound-gen}) and the fact that both $d\left(0, \partial^>\phi(\bar{x})\right)$ and $\ebm(\phi, \bar{x})$ reflect only local properties of $\phi$ near $\bar{x}$, we get the inequality $d\left(0, \partial^>\phi(\bar{x})\right)\leq \ebm(\phi, \bar{x})$ immediately. In view of (\ref{sharp-bound-convex}) and the fact that $\sigma_{\partial \phi(\bar{x})}:\R^n\rightarrow \R$ is continuous and sublinear (hence convex) as $\partial \phi(\bar{x})$ is a nonempty compact and convex set, we have $\ebm(\sigma_{\partial \phi(\bar{x})}, 0)=d\left(0, \partial^>\sigma_{\partial \phi(\bar{x})}(0)\right)$. Therefore, to show the inequality $\ebm(\phi, \bar{x}) \leq d\left(0,  \partial^>\sigma_{\partial \phi(\bar{x})}(0) \right)$, it suffices to show $\ebm(\phi, \bar{x}) \leq \ebm( \sigma_{\partial \phi(\bar{x})}, 0)$. This can be done by establishing that if there exist some $\tau>0$ and some   neighborhood $O$ of $\bar{x}$ such that
 \begin{equation}\label{sgl}
\tau d(x,[\phi\leq 0])\leq \phi(x)_+\quad\forall x\in O,
\end{equation}
then the following condition holds:
 \[
 \tau d(w, [h\leq 0])\leq h(w)_+\quad\forall w\in U,
 \]
 where $h:=\sigma_{\partial \phi(\bar{x})}$, and $U$ is a neighborhood of the origin of $\R^n$.

    Without loss of generality, we may assume that $O$ is open and that  $\phi$  is locally Lipschitz continuous on $O$. Let $U= \R^n$. In what follows, let $w\in \R^n$ be arbitrarily given with $h(w)>0$.  Since $\phi$ is regular and locally Lipschitz continuous at  $\bar{x}$, we get from   \cite[Corollary 8.19 and Exercise 9.15]{roc} that 
 \begin{equation}\label{wdlmm}
 \lim_{t\downarrow 0}\frac{\phi(\bar{x}+t w)-\phi(\bar{x})}{t}=d\phi(\bar{x})(w)=\sigma_{\partial \phi(\bar{x})}=h(w).
 \end{equation}
As $\bar{x}$ is on the boundary of the level set $[\phi\leq 0]$ and $\phi$ is locally Lipschitz continuous at $\bar{x}$, we have $\phi(\bar{x})=0$  and thus
\[
h(w)=\lim_{t\downarrow 0}\frac{\phi(\bar{x}+t w)}{t}>0,
\]
entailing that $\phi(\bar{x}+t w)>0$ for all $t>0$ sufficiently small.
By (\ref{sgl}), we have for all  $t>0$ sufficiently small,
\begin{equation}\label{eq3}
 \tau d(\bar{x}+t w, [\phi\leq 0])\leq \phi(\bar{x}+t w).
\end{equation}
Let $\kappa(x):=d(x, [\phi\leq 0])$ for all $x\in \R^n$.  Since the distance function $\kappa$ is Lipschitz continuous on $\R^n$, we get from \cite[Exercise 9.15]{roc} that
\[
 d\kappa(\bar{x})(w)=\liminf_{t\downarrow 0}\frac{d(\bar{x}+t w, [\phi\leq 0])-d(\bar{x},[\phi\leq 0])}{t}=\liminf_{t\downarrow 0}\frac{d(\bar{x}+t w, [\phi\leq 0])}{t}.
\]
In view of (\ref{wdlmm}) and (\ref{eq3}), we   have
\begin{equation}\label{wlgzq111}
\tau d\kappa(\bar{x})(w)\leq \liminf_{t\downarrow 0}\frac{\phi(\bar{x}+t w)-\phi(\bar{x})}{t}=h(w).
\end{equation}
 By  \cite[Example 8.53]{roc}, we have $d(w, T_{[\phi\leq 0]}(\bar{x}))=d\kappa(\bar{x})(w)$,
which, together with (\ref{wlgzq111}), implies that $ \tau d(w, T_{[\phi\leq 0]}(\bar{x}))\leq h(w)$. By definition, it is easy to verify that $T_{[\phi\leq 0]}(\bar{x})\subset [h\leq 0]$ and  hence $d(w, [h\leq 0])\leq d(w, T_{[\phi\leq 0]}(\bar{x}))$.
Therefore, we have $\tau d(w, [h\leq 0])\leq h(w)$.  This completes the proof. \hfill$\Box$

In view of Theorem \ref{theo-end-set-generel},  we have
\[
\partial^>\sigma_{\partial \phi(\bar{x})}(0)=\cl(\bigcup_{\sigma_{\partial \phi(\bar{x})}(w)>0}\partial \sigma_{\partial \phi(\bar{x})}(w))=\cl(\es(\partial\phi(\bar{x}))),
\]
and
\[
d(0, \partial^>\sigma_{\partial \phi(\bar{x})}(0))=d(0, \bigcup_{\sigma_{\partial \phi(\bar{x})}(w)>0}\partial \sigma_{\partial \phi(\bar{x})}(w))=d(0, \es(\partial\phi(\bar{x}))).
\]
That is, the upper estimate $d\left(0,  \partial^>\sigma_{\partial \phi(\bar{x})}(0) \right)$ in (\ref{lower-upper-esti}) is nothing else but the distance of 0 from the end set of $\partial\phi(\bar{x})$, or equivalently,  the distance of 0 from the union of all the exposed faces of $\partial\phi(\bar{x})$ having normal vectors at which the support function $\sigma_{\partial\phi(\bar{x})}$ takes positive values.

The following examples show that both the lower estimate and upper estimate in (\ref{lower-upper-esti}) may not be tight, where the first example is taken from \cite{studward99} (see also \cite{meng}) and the second one is taken from \cite[Remark 3.6]{ca2}.

\begin{example}\label{stu_war}(underestimated lower estimate).
Let $\bar{x}=0$ and let $\phi:\R\rightarrow \R_+$ be defined by
$$
\phi(x)=\left\{
\begin{array}{ll}
 0 & \mbox{if}\;x\leq 0,   \\
2^{-n}  & \mbox{if}\; 2^{-n-1}\leq x \leq 2^{-n}\;\mbox{with}\; n \;\mbox{being an odd integer},\\
3x-2^{-n}  & \mbox{if}\; 2^{-n-1}\leq x \leq 2^{-n}\;\mbox{with}\; n \;\mbox{being an even integer},\\
x  & \mbox{otherwise}.\\
\end{array}
\right.
$$
It is clear to see that $\phi$ is Lipschitz continuous and regular at $\bar{x}=0$. By some direct calculations, we have  $\partial \phi(\bar{x})=\partial^> \phi(\bar{x})=[0,1]$,  $\partial^>\sigma_{\partial \phi(\bar{x})}(0)=\es(\partial \phi(\bar{x}))=\{1\}$,  and $\ebm(\phi, \bar{x})=1$. It then follows that
\[
0=d\left(0, \partial^>\phi(\bar{x})\right)<\ebm(\phi, \bar{x}) =d\left(0,  \partial^>\sigma_{\partial \phi(\bar{x})}(0)\right)=1.
\]
That is, the lower estimate in (\ref{lower-upper-esti}) is underestimated.
\end{example}

\begin{example}(overestimated upper estimate).\label{ex-lmh}
Let $\bar{x}=(0, 0)^T$, and let $$\phi(x)=\max\{f_1(x), f_2(x)\},$$ where $f_1(x)=x_1^2+x_2^2+\frac{1}{2}(x_1+x_2)$ and $f_2(x)=x_1+x_2$. It is clear that $\phi$ is a convex function.  Clearly, $\partial \phi(\bar{x})=\co\{(\frac{1}{2},\frac{1}{2})^T, (1,1)^T\}$. From Corollary \ref{conv-hull-compa}, it follows that  $\partial^>\sigma_{\partial \phi(\bar{x})}(0)=\es(\partial\phi(\bar{x}))=\{(1, 1)^T\}$. But from Remark 3.6 (i) of \cite{ca2}, we get $\partial^>\phi(\bar{x})=\co\{(\frac{1}{2},\frac{1}{2})^T, (1,1)^T\}$. Therefore,
\[
\frac{\sqrt{2}}{2}=d(0, \partial^>\phi(\bar{x}))=\ebm(\phi, \bar{x})<d(0, \partial^>\sigma_{\partial \phi(\bar{x})}(0))=\sqrt{2}.
\]
That is, the upper estimate in (\ref{lower-upper-esti}) is overestimated.
\end{example}

\subsection{Sharp Lower Estimation for Lower $\mathcal{C}^1$ Functions}

Many functions expressed by pointwise max   of infinite collections of smooth functions
have the `subsmoothness'  property, which is between local Lipschitz continuity and strict differentiability. Our aim in this subsection is to show that the lower estimate  in (\ref{lower-upper-esti}) is a tight one for  lower $\mathcal{C}^1$ functions.

Throughout this subsection,  let  $\phi$ be lower-$\mathcal{C}^1$ on an open subset $O$ of $\R^n$ (cf. \cite[Definition 10.29]{roc}) and let $\bar{x}\in O$ be  a  fixed  point on the boundary of the level set $[\phi\leq 0]$. Moreover, we assume that on some open neighborhood $V$ of   $\bar{x}$ there is a representation
\begin{equation}\label{subsmooth}
\phi(x)=\max_{y\in Y} f(x,y)
\end{equation}
in which the functions $f(\cdot, y)$ are of class $\mathcal{C}^1$ on $V$ and
the index set $Y\subset\R^m$ is a   compact space  such that $f(x, y)$  and $\nabla_x f(x, y)$ depend continuously not just on $x\in V$ but jointly on $(x,y)\in V\times Y$. In what follows, we shall show that the lower estimate  $d(0, \partial^>\phi(\bar{x}))$ in (\ref{lower-upper-esti}) is equal to  the error bound modulus $\ebm(\phi,\bar{x})$.

To begin with, we list some nice properties of $\phi$ as follows (cf. \cite[Theorem 10.31]{roc}).
 \begin{description}
\item[(a)] $\phi$ is locally Lipschitz continuous and regular on $O$.
\item[(b)] $\partial \phi(x)=\co\{\nabla_x f(x,y)|y\in Y(x)\}$ for all $x\in V$, where  $Y:V\rightrightarrows \R^m$ is the active index set mapping defined by
\begin{equation}\label{aism}
Y(x):=\{y\in Y|f(x,y)=\phi(x)\}.
\end{equation}
\item[(c)] $\sigma_{\partial\phi(x)}(w)=d\phi(x)(w)=\displaystyle\max_{y\in Y(x)}\langle\nabla_x f(x,y),w\rangle$ for all $x\in V$ and $w\in \R^n$.
\item[(d)] The set-valued mapping $Y$ defined by (\ref{aism}) is  outer semicontinuous  at $\bar{x}$, i.e.,
\[
\limsup_{x\rightarrow \bar{x}}Y(x)\subset Y(\bar{x}).
\]
\end{description}

Next we obtain some equivalent properties for $\phi$ defined by (\ref{subsmooth}) having a local error bound.

\begin{proposition}\label{prop-equi}
Let $\tau>0$ and let
\[
\mathcal{Y}(\bar{x}):=\{Y'\subset Y(\bar{x})\mid \exists\, \{x_k\}\subset [\phi>0]\;\textup{with}\; x_k\to\bar{x}\;\textup{and}\;Y(x_k)\to Y'\}.
\]
The following properties are equivalent:
\begin{description}
\item[(i)]There exists some $\varepsilon>0$ such that for all $x\in \R^n$ with $\|x-\bar{x}\|\leq \varepsilon$,
\begin{equation}\label{eq5}
\tau d(x,[\phi\leq 0])\leq \phi(x)_+.
\end{equation}
\item[(ii)]
 For every $Y'\in \mathcal{Y}(\bar{x})$, there exists some $u\in \R^n$ with $\|u\|=1$ such that
\[
\langle\nabla_x f(\bar{x},y), u\rangle\geq \tau\quad\forall y\in Y'.
\]
\item[(iii)] For every $Y'\in \mathcal{Y}(\bar{x})$, $d(0,\co\{\nabla_x f(\bar{x},y)|y\in Y'\})\geq \tau$.
\item[(iv)] There exists some $\delta>0$ such that the inequality $d(0,\partial \phi(x))\geq \tau$  holds for all $x\in \R^n$ with $\phi(x)>0$ and  $||x-\bar{x}||\leq \delta$.
\end{description}
\end{proposition}
\noindent{\bf Proof.} For the sake of notation simplicity, we use  $C$ to denote the level set $[\phi\leq 0]$ in what follows. We shall prove  step by step that  (i)$\Longrightarrow$(ii)$\Longrightarrow$(iii)$\Longrightarrow$(iv)$\Longrightarrow$(i).

[(i)$\Longrightarrow$(ii)]: Assume that there exists some $\varepsilon>0$ such that (\ref{eq5}) holds for all $x\in \R^n$ with $\|x-\bar{x}\|\leq \varepsilon$.   First, we  show that for any $x\in \bd C\cap B(\bar{x},\frac{\varepsilon}{2})$ and any  proximal normal vector $u$ to $C$ at $x$ with $||u||=1$,  there exists some $y\in Y(x)$ such that
 \begin{equation}\label{eq4}
\langle\nabla_x f(x,y), u\rangle\geq \tau.
\end{equation}
By the definition of proximal normal vectors,  there exist some $x'\in \R^n$ and $\beta> 0$ such that
\[
u=\beta(x'-x)\;\mbox{and}\;x\in P_C(x').
\]
Take $\rho:=\min\{\frac{\varepsilon}{2}, \|x'-x\|\}$. Then it is easy to verify that
\[
x+tu\in B(\bar{x},\varepsilon)\;\forall t\in (0,\rho]\;\mbox{and}\;x\in P_C(x+tu)\;\forall t\in (0,\rho].
\]
In view of (\ref{eq5}), we have
\[
\tau t=\tau||x+tu-x||=\tau d(x+tu,C)\leq \phi(x+tu)_+\;\forall t\in (0,\rho].
\]
Thus, we have $\tau \leq \liminf_{t\to 0_+}\frac{\phi(x+tu)_+-\phi(x)_+}{t}$.
From \cite[Theorems 9.16 and 10.31]{roc},  it follows  that $\phi(x)_+$ is locally Lipschitz continuous with
\[
\liminf_{t\to 0_+}\frac{\phi(x+tu)_+-\phi(x)_+}{t}=d\phi(x)(u)_+=\max\{\max_{y\in Y(x)}\langle\nabla_x f(x,y),u\rangle,0\}.
\]
Therefore, we have $\tau \leq \max\{\max_{y\in Y(x)}\langle\nabla_x f(x,y),u\rangle,0\}$. In view of $\tau>0$, we   have
$\tau \leq \max_{y\in Y(x)}\langle\nabla_x f(x,y),u\rangle$. Since $Y(x)$ is compact, there exists some $y\in Y(x)$ such that (\ref{eq4}) holds.

Next, we show (ii) by virtue of the previous result.  Let $Y'\in \mathcal{Y}(\bar{x})$.
By definition,  there exists some sequence $\{x'_k\}\in \R^n\setminus C$ with $x'_k\to \bar{x}$ and $Y(x'_k)\to Y'$, entailing that each $y\in Y'$ corresponds to  a sequence $y'_k\to y$  such that $y'_k\in Y(x'_k)$ for all $k$. Since $C$ is a closed set,   there exists some $x_k\in \bd C$ such that $x_k\in P_C(x'_k)$.   Clearly, $x_k\rightarrow \bar{x}$  and  $u_k:=\frac{x'_k-x_k}{||x'_k-x_k||}$ is a proximal normal vector to $C$ at $x_k$.   By taking a subsequence if necessary, we can assume that $u_k\to u$, implying that $||u||=1$. In what follows, let $y\in Y'$ be given arbitrarily.  To show (ii), it suffices to show
  \begin{equation}\label{kuga}
\langle\nabla_x f(\bar{x},y),u\rangle\geq \tau.
\end{equation}
According to the previous result,  we can find some  $y_k\in Y(x_k)$ such  that for  all sufficiently large $k$,
\begin{equation}\label{eq7}
\langle\nabla_x f(x_k,y_k), u_k\rangle\geq \tau.
\end{equation}
Since all $Y(x_k)$ are subsets of the compact set $Y$, by taking a subsequence if necessary,  we can assume that $y_k\to \bar{y}$.  By the mean value theorem, there is some $\theta_k\in [0,1]$ such that
\[
f(x'_k,y'_k)-f(x_k,y'_k)=\left\langle\nabla_x f\left(x_k+\theta_k (x'_k-x_k),y'_k\right),\;x'_k-x_k\right\rangle,
\]
which, by the continuity of $\nabla_x f$,  implies that
\[
\begin{array}{l}
\;\;\displaystyle\frac{|f(x'_k,y'_k)-f(x_k,y'_k)-\langle\nabla_x f(x_k,y'_k),x'_k-x_k\rangle|}{||x'_k-x_k||}\\[0.5cm]
\displaystyle=\frac{\left\langle\nabla_x f\left(x_k+\theta_k (x'_k-x_k),y'_k\right)-\nabla_x f(x_k, y'_k),\;x'_k-x_k\right\rangle}{||x'_k-x_k||}\\[0.5cm]
\displaystyle\leq\|\nabla_x f\left(x_k+\theta_k (x'_k-x_k),y'_k\right)-\nabla_x f(x_k, y'_k)\| \to 0.
\end{array}
\]
Thus, we have
\begin{equation}\label{fangxin}
\lim_{k\to +\infty}\frac{f(x'_k,y'_k)-f(x_k,y'_k)}{||x'_k-x_k||}=\lim_{k\to +\infty}\langle\nabla_x f(x_k,y'_k), u_k\rangle=\langle\nabla_x f(\bar{x},y),u\rangle.
\end{equation}
Similarly, we  obtain
\begin{equation}\label{yefangxin}
\lim_{k\to +\infty}\frac{f(x'_k,y_k)-f(x_k,y_k)}{||x'_k-x_k||}=\lim_{k\to +\infty}\langle\nabla_x f(x_k,y_k), u_k\rangle\geq \tau,
\end{equation}
where the inequality follows from (\ref{eq7}).  Observing that
\[
f(x'_k,y'_k)-f(x_k,y'_k)\geq \phi(x'_k)-\phi(x_k)\geq f(x'_k,y_k)-f(x_k,y_k),
\]
we get from (\ref{fangxin}) and (\ref{yefangxin}) that (\ref{kuga}) holds. This completes the proof for (i)$\Longrightarrow$(ii).

[(ii)$\Longrightarrow$(iii)]:  Let $Y'\in \mathcal{Y}(\bar{x})$. By (ii), there exists  some $u\in \R^n$ with $||u||=1$ that
\[
\langle u, v\rangle\geq \tau\geq \langle u, w\rangle,\;\forall v\in \co\{\nabla_x f(\bar{x},y)|y\in Y'\}, \forall w\in  B(0,\tau).
\]
Then by a separation argument, we have
\[
0\not\in \inte(\co\{\nabla_x f(\bar{x},y)|y\in Y'\}-B(0,\tau)),
\]
which clearly  implies  (iii).

[(iii)$\Longrightarrow$(iv)]: Let $\tau'\in (0, \tau)$ be given arbitrarily. First, we shall prove  that,  there exists some $\delta>0$ such that for all $x\not\in C$ with $||x-\bar{x}||\leq \delta$,
\begin{equation}\label{eq13}
d(0,\co\{\nabla_x f(x,y)|y\in Y(x)\})\geq \tau'.
\end{equation}
Suppose by contradiction that (\ref{eq13}) does not hold, i.e., there exists a sequence $\{x_k\}\subset \R^n\setminus C$ with $x_k\to \bar{x}$ and
\[
d(0,\co\{\nabla_x f(x_k,y)|y\in Y(x_k)\})< \tau'.
\]
It follows from the Carath\'{e}odory theorem that,   there exist some $t_k^j\geq 0$ and $y_k^j\in Y(x_k)$ with  $j=1,2,\cdot \cdot \cdot,n+1$  such that
\begin{equation}\label{eq14}
\Sigma_{j=1}^{n+1}t_k^j=1\;\mbox{and}\; ||\Sigma_{j=1}^{n+1}t_k^j\nabla_x f(x_k,y_k^j)||\leq \tau'.
\end{equation}
Since $Y(x_k)\subset Y$ for all $k$ and $Y$ is compact, it follows from \cite[Theorem 4.18]{roc} that $Y(x_k)$ has a subsequence converging to $Y^*$, a subset of $Y$. By taking a subsequence if necessary, we assume that
\[
Y(x_k)\to Y^*,\quad t_k^j\to t^j\geq 0,\quad\mbox{and}\quad y_k^j\to y^j\in Y^*.
\]
Since $Y: V \rightrightarrows \R^m$ defined by (\ref{aism}) is osc at $\bar{x}$, it follows from \cite[Exercise 5.3]{roc} that $Y^*\subset Y(\bar{x})$, entailing that  $Y^*\in \mathcal{Y}(\bar{x})$.
By (\ref{eq14}) and the continuity of $\nabla_x f$, we  have
\[
\Sigma_{j=1}^{n+1}t^j=1\;\mbox{and}\; ||\Sigma_{j=1}^{n+1}t^j\nabla_x f(\bar{x},y^j)||\leq \tau'.
\]
Thus, we have  $d(0,\co\{\nabla_x f(\bar{x},y)|y\in Y^*\})\leq \tau'$, contradicting  to (ii). This contradiction implies that (\ref{eq13}) holds. Since $\tau'\in (0, \tau)$  is given arbitrarily,  we  confirm that there exists some $\delta>0$ such that the following inequality  holds for all $x\not\in C$ with $||x-\bar{x}||\leq \delta$:
\begin{equation}\label{eq12}
d(0,\co\{\nabla_x f(x,y)|y\in Y(x)\})\geq \tau.
\end{equation}
    In view of (b),  we can reformulate (\ref{eq12})  as $d(0,\partial \phi(x))\geq \tau$.

[(iv)$\Longrightarrow$(i)]: This implication follows readily from  \cite[Proposition 2.1]{meng}.    \hfill $\Box$

\begin{remark}
 When $Y$ is a finite set,  the results in Proposition \ref{prop-equi} can be found in \cite[Theorem 2.1]{meng}. See also Kummer \cite{ku}. In the semi-infinite setting, Proposition \ref{prop-equi} improves the corresponding results in Henrion and Outrata \cite{hen} and Zheng and Yang \cite{zheng2}.
 \end{remark}

Next theorem shows that the lower estimate in (\ref{lower-upper-esti}) is a tight one.

\begin{theorem}\label{theo-sip}
The following equalities hold:
\begin{equation}\label{out-lim-subsmooth}
\partial^>\phi(\bar{x})=\bigcup_{Y'\in \mathcal{Y}(\bar{x})}\co\{\nabla_x f(\bar{x},y)|y\in Y'\},
\end{equation}
and
\[
{\rm ebm}(\phi, \bar{x})=d(0, \partial^>\phi(\bar{x})).
\]

\end{theorem}
\noindent{\bf Proof.} The equality (\ref{out-lim-subsmooth}) follows readily from the definition of outer limiting subdifferential and the fact that all $Y(x)$ are compact and convex subsets of $Y(\bar{x})$ when $x$ is close enough to $\bar{x}$. The equality ${\rm ebm}(\phi, \bar{x})=d(0, \partial^>\phi(\bar{x}))$ follows from (\ref{out-lim-subsmooth}) and the equivalence of (i) and (iii) in Proposition \ref{prop-equi}.  \hfill $\Box$

The upper estimate  $d(0, \partial^>\sigma_{\partial \phi(\bar{x})}(0))$ in (\ref{lower-upper-esti}) has an alternative expression  in terms of a collection of subsets of the index set $Y(\bar{x})$ defined by \[
\mathcal{Y}^>(\bar{x}):=\left\{\,Y'\subset Y(\bar{x})\left| \, \exists w\in \R^n: Y'=\argmax_{y\in Y(\bar{x})}\left\langle \nabla_x f(\bar{x},y), w\right\rangle,\,\max_{y\in Y(\bar{x})}\left\langle \nabla_x f(\bar{x},y), w\right\rangle>0\right.\right\}.
\]
By applying Corollary \ref{conv-hull-compa}, we have
\begin{equation}\label{tenggeer111}
\begin{array}{l}
\displaystyle\bigcup_{Y'\in \mathcal{Y}^>(\bar{x})}  \co\{\nabla_x f(\bar{x},y)|y\in Y'\}\subset \es(\partial \phi(\bar{x}))=\gamma_{\partial \phi(\bar{x})}^{-1}(1)\\
\;\;\subset\cl\left(\displaystyle\bigcup_{Y'\in \mathcal{Y}^>(\bar{x})}  \co\{\nabla_x f(\bar{x},y)|y\in Y'\}\right)=\partial^>\sigma_{\partial \phi(\bar{x})}(0),
\end{array}
\end{equation}
  where each $\co\{\nabla_x f(\bar{x},y)|y\in Y'\}$   is an exposed face of $\partial \phi(\bar{x})$. Thus
  \[d(0, \partial^>\sigma_{\phi(\bar{x})}(0))=d(0,\bigcup_{Y'\in \mathcal{Y}^>(\bar{x})}\co\{\nabla_x f(\bar{x},y)|y\in Y'\}).
  \]
  If the index set $Y(\bar{x})$ is finite, all the inclusions in (\ref{tenggeer111}) become equalities.

\subsection{Sharp Upper Estimation for Convex Functions}

In the case of $\phi$ being  finite and convex on some convex neighborhood of $\bar{x}$, entailing that  $\phi$ is regular and locally Lipschitz continuous on some open neighborhood of $\bar{x}$  (cf. \cite[Examples 7.27 and 9.14]{roc}), the lower estimate in (\ref{lower-upper-esti}) is tight,  but the upper estimate in (\ref{lower-upper-esti}) could be overestimated,  as seen in Example \ref{ex-lmh}.

In general, we cannot expect that the upper estimate in (\ref{lower-upper-esti})  is a tight one,   unless some regularity conditions are imposed as we have done in the following theorem.

\begin{theorem} \label{theo-sharp-upper}
 Assume that $\phi$ is finite and convex on some convex neighborhood of $\bar{x}$. If there is a neighborhood $V$ of $\bar{x}$ such that
 \begin{equation}\label{liangtou}
  [\phi\leq 0]\cap V=(\bar{x}+[d\phi(\bar{x})\leq 0])\cap V,
  \end{equation}
   then the following equalities hold:
 \begin{equation}\label{both-sharp-conv}
d\left(0, \partial^>\phi(\bar{x})\right)=\ebm(\phi, \bar{x}) = d\left(0,  \partial^>\sigma_{\partial \phi(\bar{x})}(0) \right).
\end{equation}
\end{theorem}
\noindent{\bf Proof.} By  the definitions of tangent cone and subderivative, we can easily verify that $[\phi\leq 0]\subset \bar{x}+T_{[\phi\leq 0]}(\bar{x})$ and $T_{[\phi\leq 0]}(\bar{x}) \subset [d\phi(\bar{x})\leq 0]$. Thus,  the regularity condition (\ref{liangtou}) amounts to that
\begin{equation}\label{acq-property}
[d\phi(\bar{x})\leq 0]=T_{[\phi\leq 0]}(\bar{x}),
\end{equation}
 and
\begin{equation}\label{eta-property}
[\phi\leq 0]\cap V=(\bar{x}+T_{[\phi\leq 0]}(\bar{x}))\cap V.
\end{equation}
 In view of (\ref{sharp-bound-convex}) and the assumption that $\phi$ is finite and convex on some convex neighborhood of $\bar{x}$, we get the first equality in (\ref{both-sharp-conv}) immediately.  It remains to show that the second equality holds under (\ref{acq-property}) and (\ref{eta-property}).  Without loss of generality, we   assume that there exists an open ball $O:=\{x\in \R^n\mid \|x\|<\delta\}$ of radius $\delta>0$ such that $\phi$ is finite and convex on $\bar{x}+O$ and that $V$ in (\ref{eta-property}) can be replaced by $\bar{x}+O$. As $\phi$ is assumed to be  finite and convex on some convex neighborhood of $\bar{x}$, it follows from \cite[Examples 7.27 and 9.14, Theorem 9.16]{roc} that  $\phi$ is regular and locally Lipschitz continuous on some open neighborhood of $\bar{x}$, and hence that $\sigma_{\partial \phi(\bar{x})}=d\phi(\bar{x})$ and
\begin{equation}\label{con-lip}
\ebm(d\phi(\bar{x}), 0)=\ebm(\sigma_{\partial \phi(\bar{x})}, 0)=d\left(0,  \partial^>\sigma_{\partial \phi(\bar{x})}(0) \right).
 \end{equation}
 Moreover, we get from Theorem \ref{theo-upper-bound}  that  $\ebm(\phi, \bar{x})\leq d\left(0,  \partial^>\sigma_{\partial \phi(\bar{x})}(0) \right)$. In the case of $d\left(0,  \partial^>\sigma_{\partial \phi(\bar{x})}(0) \right)=0$,   the second equality in (\ref{both-sharp-conv}) holds trivially. So in what follows we assume that $d\left(0,  \partial^>\sigma_{\partial \phi(\bar{x})}(0) \right)>0$.

Let $0<\tau<d\left(0,  \partial^>\sigma_{\partial \phi(\bar{x})}(0) \right)$. In view of (\ref{con-lip}) and the positive homogeneity of $d\phi(\bar{x})$, the following condition holds:
\begin{equation}\label{xianggexiaohai}
\tau d(w, [d\phi(\bar{x})\leq 0])\leq d\phi(\bar{x})(w)_+\;\forall w\in \R^n.
 \end{equation}
Let $x\in\bar{x}+\frac{1}{2} O$ be arbitrarily chosen. It is straightforward to verify that
\[
d(x-\bar{x}, T_{[\phi\leq 0]}(\bar{x}))= d(x-\bar{x}, T_{[\phi\leq 0]}(\bar{x})\cap O)=d(x,(\bar{x}+T_{[\phi\leq 0]}(\bar{x}))\cap (\bar{x}+O)),
\]
and
\[
d(x,[\phi\leq 0])=d(x,[\phi\leq 0]\cap (\bar{x}+O)).
\]
In view of  (\ref{eta-property}), we have
\[
d(x,[\phi\leq 0])=d(x-\bar{x}, T_{[\phi\leq 0]}(\bar{x})),
\]
which implies by   (\ref{acq-property}) that
\[
d(x,[\phi\leq 0])\leq d(x-\bar{x}, [d\phi(\bar{x})\leq 0]).
\]
By (\ref{xianggexiaohai}), we have
\begin{equation}\label{guoer}
\tau d(x,[\phi\leq 0])\leq d\phi(\bar{x})(x-\bar{x})_+.
\end{equation}
Since $\phi$ is finite and convex on $\bar{x}+O$, we get from
\cite[Proposition 8.21]{roc} that
\begin{equation}\label{guoer111}
d\phi(\bar{x})(x-\bar{x})\leq \phi(x)-\phi(\bar{x})=\phi(x).
\end{equation}
 In view of
(\ref{guoer}) and (\ref{guoer111}), we have $\tau d(x,[\phi\leq 0])\leq \phi(x)_+$.
Since $x\in\bar{x}+\frac{1}{2} O$ is chosen arbitrarily, we thus have $ \tau\leq \ebm(\phi, \bar{x})$, entailing that $d\left(0,  \partial^>\sigma_{\partial \phi(\bar{x})}(0) \right)\leq \ebm(\phi, \bar{x})$.   This completes the proof. \hfill$\Box$

\begin{remark}\label{rem-acq+eta}
 Recall that the Abadie constraint qualification \cite{li} (ACQ, for short) holds at $\bar{x}$ if (\ref{acq-property}) holds, and that the level set $[\phi\leq 0]$ admits exactness of tangent approximation (ETA, for short) at  $\bar{x}$ if there exists some neighborhood $V$ of $\bar{x}$ such that (\ref{eta-property}) holds. From the proof of Theorem \ref{theo-sharp-upper}, it is clear that the regularity condition (\ref{liangtou})
amounts to  the ACQ  plus the ETA.  It turns out  in last section that,  the outer limiting subdifferential set $\partial^>\sigma_{\partial \phi(\bar{x})}(0)$, unlike the outer limiting subdifferential set $\partial^>\phi(\bar{x})$, depends on the nominal point $\bar{x}$ only and does not get the nearby points involved. As can be seen from Theorem  \ref{theo-sharp-upper}, it is  the ETA property that makes it possible for  $d(0, \partial^>\sigma_{\partial \phi(\bar{x})}(0))$ to serve as the error bound modulus $\ebm(\phi, \bar{x})$ which normally depends on not only $\bar{x}$ but its nearby points. Note that the idea of using the ETA property has already appeared in Zheng and Ng \cite{zheng1} and that various characterizations of the ETA property has been presented in \cite{meng1}.  If the ETA property (\ref{eta-property}) does not hold,  the upper estimate $d(0, \partial^>\sigma_{\partial \phi(\bar{x})}(0))$ may be overestimated as can be seen from Example \ref{ex-lmh}, in which $[\phi\leq 0]=\{x\in \R^2| x_1^2+x_2^2+\frac{1}{2}(x_1+x_2)\leq 0\}$ and the ETA property doest not hold at any $x\in [\phi\leq 0]$.
 \end{remark}

In the remainder of this subsection, we apply   Theorem \ref{theo-sharp-upper} to  the linear system
\begin{equation}\label{lop-sys}
\langle a_t, x\rangle\leq b_t\quad\forall t\in T,
\end{equation}
where $a_t\in \R^n$, $b_t\in \R$, and $T$ is  a  compact  space such that $a_t$ and $b_t$ depend continuously on $t\in T$. In what follows, let $\phi(x):=\max_{t\in T}\{\langle a_t, x\rangle -b_t\}$ and let
 $T(x):=\{t\in T\mid \langle a_t, x\rangle -b_t=\phi(x)\}$.  Clearly, the level set $[\phi\leq 0]$ is the solution set of the linear system (\ref{lop-sys}), and the regularity condition (\ref{liangtou}) specified for   $x\in [\phi\leq 0]$  can be reformulated as
   \begin{equation}\label{liangtou111}
   \{y\mid \langle a_t, y\rangle\leq b_t\;\forall t\in T\}\cap V=(x+\{w\mid \langle a_t, w\rangle\leq 0\;\forall t\in T(x)\})\cap V,
  \end{equation}
  where $V$ is a neighborhood of $x$.

 Our first result for the linear system (\ref{lop-sys}) assumes the  regularity  condition (\ref{liangtou111}) on one nominal point in the solution set only.

 \begin{corollary}\label{coro-LOP}
Consider a solution $x$ to the linear system (\ref{lop-sys}). If the regularity condition (\ref{liangtou111}) holds, then
\begin{equation}\label{lop-modulus}
d\left(0, \partial^>\phi(x)\right)=\ebm(\phi, x) = d\left(0,  \partial^>\sigma_{\partial \phi(x)}(0) \right)=d(0, \bigcup_{T'\in \mathcal{T}(x)}\co\{a_t\mid t\in T'\}),
\end{equation}
where
\[
\mathcal{T}(x):=\{T'\subset T(x)\mid \exists w\in \R^n: \langle a_t, w\rangle =1\,\forall t\in T',\;\langle a_t, w\rangle<1\,\forall t\in T(x)\backslash T'\}.
\]
\end{corollary}
\noindent{\bf Proof.} Applying  Theorem \ref{theo-sharp-upper}, we get the first two equalities in (\ref{lop-modulus}). Applying  Corollary \ref{conv-hull-compa}, we get the third equality in (\ref{lop-modulus}) by taking Remark \ref{rem-canovas} into account.  This completes the proof. \hfill$\Box$

Our second result for  the linear system (\ref{lop-sys})  assumes the regularity condition (\ref{liangtou111}) on the whole solution set, leading to a locally polyhedral linear system  as defined in \cite{agl98}, which requires that
  \begin{equation}\label{lop-con}
\left(\pos\co\{a_t\mid t\in T(x)\}\right)^*=\pos([\phi\leq 0]-x)\quad \forall x\in [\phi\leq 0].
 \end{equation}
  As a finite linear system is naturally  locally polyhedral, our result below recovers \cite[Theorem 4.1]{ca} for the case of a finite linear system.
 \begin{corollary}\label{coro-LOP}
Consider the linear system (\ref{lop-sys}). The equalities in (\ref{lop-modulus}) hold for all $x\in [\phi\leq 0]$ if  one of the following equivalent properties is satisfied:
\begin{description}
  \item[(a)] The regularity condition (\ref{liangtou111}) holds for all $x$ in the solution set $[\phi\leq 0]$;
  \item[(b)] The linear system (\ref{lop-sys}) is locally polyhedral, i.e., (\ref{lop-con}) holds.
\end{description}
\end{corollary}
\noindent{\bf Proof.} It suffices to show the equivalence of (a) and (b). To begin with, we point out that   $d\phi(x)(w)=\max_{t\in T(x)}\langle a_t, w\rangle$ as can be seen from \cite[Theorem 10.31]{roc},   and that $[\phi\leq 0]$ is convex (implying that $T_{[\phi\leq 0]}(x)=\cl \pos([\phi\leq 0]-x)$).  Moreover, we have
 \begin{equation}\label{gzjj000}
 [\phi\leq 0]-x\subset \pos([\phi\leq 0]-x)\subset  T_{[\phi\leq 0]}(x)\subset [d\phi(x)\leq 0],
\end{equation}
 and

  \begin{equation}\label{gzjj}
  \begin{array}{ll}
  \left(\pos\co\{a_t\mid t\in T(x)\}\right)^* &=\{a_t\mid t\in T(x)\}^* \\
&=\{w\in \R^n\mid \langle a_t, w\rangle\leq 0\;\forall t\in T(x)\}\\
&=[d\phi(x)\leq 0].
  \end{array}
  \end{equation}

 First, we show  $(b)\Longrightarrow (a)$.
    Condition (\ref{lop-con}) implies that $\pos([\phi\leq 0]-x)$ is closed for all $x\in [\phi\leq 0]$.
  In view of \cite[Proposition 4.1]{meng1}, the level set $[\phi\leq 0]$ admits the ETA property (\ref{eta-property}) at every $x\in [\phi\leq 0]$.    By (\ref{lop-con}) and  (\ref{gzjj}),  the ACQ (\ref{acq-property}) holds for all $x\in [\phi\leq 0]$. In view of Remark \ref{rem-acq+eta}, the regularity condition (\ref{liangtou}) or its reformulation (\ref{liangtou111}) holds for all $x\in [\phi\leq 0]$.

Now we show $(a)\Longrightarrow (b)$. Let $x\in [\phi\leq 0]$. Assume that the regularity condition (\ref{liangtou111}) or its reformulation (\ref{liangtou}) holds. It then follows from (\ref{gzjj000}) that
\[
\pos([\phi\leq 0]-x)=[d\phi(x)\leq 0],
\]
which together with (\ref{gzjj}) implies (\ref{lop-con}). This completes the proof. \hfill$\Box$

To end this subsection, we illustrate two examples selected from  \cite{ca}.  By Example \ref{not-lop}, we demonstrate that (\ref{lop-modulus}) may not hold  if the linear system (\ref{lop-sys}) is not locally polyhedral, and by Example \ref{lop-not}, we demonstrate that (\ref{lop-modulus}) may still hold even if  the linear system (\ref{lop-sys}) is not locally polyhedral.

\begin{example}\label{not-lop}
Let $\bar{x}=(1,0)^T$ and  $\phi(x)=\max_{t\in T}\{\langle a_t, x\rangle -b_t\}$, where $T=[0, 2\pi]$, $a_t=(t\cos t, t\sin t)^T$ and $b_t=t$. Clearly, $[\phi\leq 0]=\{x\in \R^n\mid \|x\|\leq 1\}$. Thus,  $\pos([\phi\leq 0]-\bar{x})$ is not closed, implying that (\ref{lop-con}) does not hold at $\bar{x}$ and that   the linear system (\ref{lop-sys}) cannot be locally polyhedral. From Example 1 of \cite{ca}, it follows that $d\left(0, \partial^>\phi(\bar{x})\right)=\ebm(\phi, \bar{x})=0$. Observing that $T(\bar{x})=\{0, 2\pi\}$ and $\mathcal{T}(\bar{x})=\{\{2\pi\}\}$, we get
\[
d\left(0,  \partial^>\sigma_{\partial \phi(\bar{x})}(0) \right)=d(0, \bigcup_{T'\in \mathcal{T}(\bar{x})}\co\{a_t\mid t\in T'\})=2\pi.
\]
That is, the upper estimate $d\left(0,  \partial^>\sigma_{\partial \phi(\bar{x})}(0) \right)$ is overestimated.
\end{example}

\begin{example}\label{lop-not}
Let $\bar{x}=(1,0)^T$ and  $\phi(x)=\max_{t\in T}\{\langle a_t, x\rangle -b_t\}$, where $T=[0, 2\pi]$, $a_t=(\cos t, \sin t)^T$ and $b_t=1$. Clearly, $[\phi\leq 0]=\{x\in \R^n\mid \|x\|\leq 1\}$. Thus,  $\pos([\phi\leq 0]-\bar{x})$ is not closed, implying that (\ref{lop-con}) does not hold at $\bar{x}$ and that   the linear system (\ref{lop-sys}) cannot be locally polyhedral. By some direct calculations, we have  $T(\bar{x})=\{0, 2\pi\}$,  $\mathcal{T}(\bar{x})=\{\{0, 2\pi\}\}$, and
\[
d\left(0,  \partial^>\sigma_{\partial \phi(\bar{x})}(0) \right)=d(0, \bigcup_{T'\in \mathcal{T}(\bar{x})}\co\{a_t\mid t\in T'\})=1.
\]
Moreover, we have $\partial\phi(\bar{x})=\co\{a_t\mid t\in T(\bar{x})\}=(1, 0)^T$ and hence $\partial^>\phi(\bar{x})=(1, 0)^T$, entailing that
\[
d\left(0, \partial^>\phi(\bar{x})\right)=\ebm(\phi, \bar{x})=1.
\]
That is,  (\ref{lop-modulus}) still holds even when the linear system (\ref{lop-sys}) is not locally polyhedral.
\end{example}

\section{Conclusions and Perspectives}

When $\phi$ is regular and locally Lipschitz continuous on some neighborhood  of $\bar{x}\in \bd([\phi\leq 0])$, we obtained in Theorem \ref{theo-upper-bound}  a lower estimate  and an upper estimate of the local error bound modulus $\ebm(\phi, \bar{x})$ as follows:
\[
d\left(0, \partial^>\phi(\bar{x})\right)\leq \ebm(\phi, \bar{x}) \leq d\left(0,  \partial^>\sigma_{\partial \phi(\bar{x})}(0) \right).
\]
In particular, when $\phi$ is finite and convex on some convex neighborhood  of $\bar{x}\in \bd([\phi\leq 0])$, we obtained in Theorem \ref{theo-sharp-upper} under the ACQ and ETA properties the following:
\[
d\left(0, \partial^>\phi(\bar{x})\right)= \ebm(\phi, \bar{x}) =d\left(0,  \partial^>\sigma_{\partial \phi(\bar{x})}(0) \right),
\]
and when $\phi$ is a lower $\mathcal{C}^1$ functions, we obtained in Theorem \ref{theo-sip} the following:
 \[
d\left(0, \partial^>\phi(\bar{x})\right)= \ebm(\phi, \bar{x})\leq d\left(0,  \partial^>\sigma_{\partial \phi(\bar{x})}(0) \right).
\]
One open question  is whether the  inclusion
\begin{equation}\label{two-outer-limiting}
\partial^>\sigma_{\partial \phi(\bar{x})}(0) \subset \partial^>\phi(\bar{x})
\end{equation}
holds or not in the general case or in some particular settings. By trying to find answers to this open question, one may need to look into the differential structure of the functions in question and need to apply some delicate modern variational tools.
It is worth noting that \cite[Theorem 3.1]{ca2} shows that  (\ref{two-outer-limiting})  holds as an equality when $\phi$ is the pointwise max of a finite collection of affine functions.
When $\phi$ is the pointwise max of a finite collection of smooth functions, \cite[Theorem 3.2]{ca2} shows that a subset of the set $\partial^>\sigma_{\partial \phi(\bar{x})}(0)$ is included in $\partial^>\phi(\bar{x})$.

\end{document}